\documentclass[a4paper, 12pt]{amsart}

\usepackage{amsmath, amssymb, amsthm}
\usepackage[english]{babel}
\usepackage[T1]{fontenc}
\usepackage{enumerate}
\def\ff{\widehat f}

\def\gg{\widehat g}

\newcommand{\ee}{\mathbf{e}}

\def\CI {\mathcal{I}}

\def \be {  \varpi}

\newcommand{\fer}[1]{(\ref{#1})}

\def\bx{{\bf x}}
\def\bX{{\bf X}}
\def\bY{{\bf Y}}

\def\by{{\bf y}}

\def\bxi{{\pmb\xi}}

\newcommand{\R}{\mathbb R}
\def\be#1\ee{\begin{equation}#1\end{equation}}

\numberwithin{equation}{section}

\newcommand{\bq}{\begin{equation}}
\newcommand{\eq}{\end{equation}}

\theoremstyle{definition}

\newenvironment{equations}{\equation\aligned}{\endaligned\endequation}

\begin{document}

\title[The Energy distance in evolution problems]{
 ENERGY DISTANCE  AND EVOLUTION PROBLEMS: \\ A PROMISING TOOL FOR KINETIC EQUATIONS}
 
 \author{GENNARO AURICCHIO}

\address{Department of Mathematics, University of Padua\\
via Trieste 36,
Padova, 35100 Italy}
\email{gennaro.auricchio@unipd.it}

\author{GIUSEPPE TOSCANI}
\address{Department of Mathematics, University of Pavia, and IMATI of the National Council for Research;
via Ferrata 1,
Pavia, 27100 Italy}
\email{giuseppe.toscani@unipv.it}

\begin{abstract}
We study the rate of convergence to equilibrium of the solutions to  Fokker--Planck type equations  with linear drift by means of Cram\'er and Energy distances,  which have been recently widely used in problems related to AI, in particular for tasks related to machine learning. In all cases in which the Fokker--Planck type equations can be treated through these distances, it is shown that the rate of decay is improved with respect to known results which are based on the decay of relative entropy.

\end{abstract}


\maketitle
\section{Introduction}\label{intro}

In recent years, Cram\'er's distance \cite{Cra}, also known as the continuous ranked probability score, has proven to be a useful tool in problems related to AI, and in particular to machine learning \cite{Belle}, especially for tasks involving distribution estimation, regression, and image generation. As explained in \cite{Belle}, the main reason behind the rediscovery of Cram\'er's distance, which is a metric for comparing probability distributions, lies in its properties, in particular that of having unbiased sample gradients, making it easier to optimize using stochastic gradient descent. 

Unfortunately, Cramér's distance can only be successfully applied to problems involving univariate probability distributions. A closely related distance, which shares the same properties as Cramér's distance \cite{ABGT2}, but is designed to also compare multivariate probability distributions, is the Energy distance introduced by Székely \cite{Sze89, Sze03, SR}.  
In \cite{ABGT2}, the relationships between the Energy distance and other well-known metrics—such as the Wasserstein distance \cite{villani2009optimal} and the Fourier-based metrics introduced in \cite{GTW}, which are widely used in kinetic theory—were thoroughly analyzed. These metrics are of particular interest due to their applicability in the field of artificial intelligence \cite{Au1,Au2,cuturi2014fast,scagliotti2025normalizing}.

In this paper, we show that these distances are particularly adapted to study the rate of convergence to equilibrium of the solutions to the linear diffusion equation and to Fokker--Planck type equations characterized by a linear drift, including the one-dimensional Fokker--Planck type equations  modeling socio-economic phenomena studied in \cite{FPTT}, providing in all cases a decay better than exponential. 

A further advantage of this analysis is related to the numerical approximation of evolution equations of Fokker--Planck type. In reason of the simple form of the Energy distance, which is essentially based on the computation of moments, the numerical study of the evolution of solutions towards a steady state could be greatly simplified.

In more details, in Section \ref{sec:basic} we briefly report the definitions and the main properties of both Cram\'er and Energy distances, as well as their connections. The equivalence between these two distances in the one-dimensional case will be shown through a  proof which is closely related to an analogous one that can be found in the literature  devoted to the Gini index \cite{Xu}.  Then, in Section \ref{sec:FP} we will show how these distances can be fruitfully employed to study convergence to equilibrium of solutions to one-dimensional Fokker--Planck equations with linear drift and general diffusion term, and to recover exponential rate of convergence in time. In Section \ref{sec:RN} we will show that the Energy distance can be used to recover exponential convergence to equilibrium in the classical linear Fokker--Planck equation in any dimension. Last, Section \ref{sec:heat} will be devoted to an application of the Energy distance to the linear diffusion equation. There, the rate of decay of the energy distance between two solutions of the linear diffusion  is shown to improve with the increasing of the value of the order $\alpha$ characterizing the distance.  Some conclusions will follow in Section \ref{sec:conclusions}.

\section{Basic facts about Cram\'er and Energy distances}\label{sec:basic}

Let $F$ be a probability measure over $\R$. When $F$ is continuous, we will assume it has density $f$.
We will associate $F$ to a random variable $X$, such that for a subset of the reals $A \subseteq \R$, we have 
\[
F(X\in A) = \int_A dF(x).
\]
The (cumulative) distribution function of $X$ is then
\be\label{cumu}
F(x) = P(X \le x) = \int_{-\infty}^x dF(x).
\ee
Given $s>0$, we denote by $P_s(\R)$ the class of all probability measures $F$ over $\R$ such that
\[
m_s(F)  = \int_{\R} |x|^s dF(x) < + \infty.
\]
If $s=1$,  the quantity $m_1(F)$ defines the (absolute) mean value of $X$, say $E(|X|)$. In this relevant case we can express the standard mean value $E(X)$ directly in terms of the cumulative function, through the formula
\be\label{for-mean}
E(X) = -\int_{-\infty}^0 F(x)\,dx + \int_0^{+\infty} (1-F(x))\, dx.
\ee
Given two random variables $X$ and $Y$ of cumulative distribution functions $F$ (respectively $G$) the Cram\'er distance between $X$ and $Y$ is given by \cite{Cra}
\be\label{Cramer}
d(F,G) = \int_\R[F(x) -G(x)]^2\, dx.
\ee 
Note that, as written, Cram\'er distance is not a metric, while its square root is. 

One of the interesting properties of Cram\'er's distance is that it can be written equivalently in different ways. It is immediate to show that it can be written in terms of Fourier transforms, where, as usual, the Fourier transform of a probability density $f$ in $P_s(\R)$ is defined by
\[
\ff(\xi) = \int_{\R}  e^{-i\xi x} \, dF(x).
\]
Indeed, for any given pair of probability measures $F,G \in P_s(\R)$,  Parseval formula implies
\be\label{Par}
\int_\R [F(x) - G(x)]^2 \, dx = \frac 1{2\pi} \int_\R |\widehat F(\xi) - \widehat G(\xi)|^2 \, d\xi,
\ee
where $\widehat F$ and $ \widehat G$ are the Fourier transforms of the probability measures $F,G$. On the other hand, if
 \[
\ff(\xi) = \int_{\R}  e^{-i\xi x} \, dF(x), \quad \gg(\xi) = \int_{\R}  e^{-i\xi x} \, dG(x),
\] 
 integration by parts gives 
 \[
 \widehat F(\xi) - \widehat G(\xi) = \frac{\ff(\xi) -\gg(\xi)}{i\xi}.
 \]
  Consequently, Cram\'er distance \fer{Cramer} can be equivalently written as
 \be\label{Cramer2}
d(F,G)  = \frac 1{2\pi} \int_\R \frac{|\ff(\xi) -\gg(\xi)|^2}{\xi^2} \, d\xi.
\ee
Also, assuming that the random variables $X$ and $Y$ are mutually independent, and that $X'$ and $Y'$ are independent copies of $X$ (respectively $Y$), Cram\'er's distance can be expressed as
\be\label{Cramer3}
d(F,G) = E(|X-Y|) - \frac 12 E(|X-X'|) - \frac 12 E(|Y-Y'|).
\ee
Note that \fer{Cramer3} implies that the Cram\'er distance of $F$ and $G$ is bounded as soon as $F,G \in P_1(\R)$.

This unexpected way of writing Cram\'er's distance, which has been shown in \cite{SR}, can be more easily derived by resorting to  an older result related to Gini inequality index \cite{Gini1,Gini2}. 
For a random variable $X \ge 0$ distributed according to $F$, and given two independent copies $X$ and $X'$, the quantity
\be\label{Gini} 
G(X)  = \frac 1{2E(X)} \int_{\R_+\times\R_+} |x-y| dF(x)dF(y) = \frac 1{2E(X)}E(|X-X'|) \le 1
\ee
defines the well-known Gini index.   Gini index can be expressed in different ways \cite{Xu}. Following the approach of Kendall and Stuart \cite{Kend},  Dorfman \cite{Dorf}  proposed  for the Gini index of a continuous income distribution the equivalent expression
\be\label{Gini2}
G(X) = 1 - \frac 1{2E(X)} \int_0^{+\infty}(1-F(x))^2\, dx.
\ee
A proof of this equivalence, that was omitted in \cite{Dorf} and in a subsequent paper on the same subject \cite{Gast}, can be found in \cite{Xu}, and relies on a simple remark, namely that, for given real numbers $x$ and $y$
\[
|x-y| = x + y -2 \min\{x,y\},
\] 
so that
\[
E(|X-Y|) = E(X)+E(Y) -2 E(\min\{X, Y\}).
\]
It is immediate to show that, for a given pair of independent random variables $X$ and $Y$,  the cumulative function $H(x)$ of $\min\{X, Y\}$ takes the form
\be\label{min}
H(x) = P[ \min\{X, Y\} \le x] = 1 - P[(X>x)\cap (Y>x)] = 1- [1-F(x)][1-G(x)].
\ee
Hence, owing to \fer{min} and to formula \fer{for-mean} we get
\begin{equations}\nonumber
E(|X-Y|) &= E(X)+E(Y)  + 2 \int_{-\infty}^0 H(x)\,dx - 2\int_0^{+\infty} (1-H(x))\, dx \\
&= E(X)+E(Y) + 2\int_{-\infty}^0[ F(x) +G(x) -F(x)G(x)]\, dx +\\
&-2\int_0^{+\infty}[ 1 -F(x)  - G(x) + F(x)G(x)]\, dx.
\end{equations}
Likewise, given two independent copies $X$ and $X'$ distributed according to $F$ we obtain
\begin{equation}\nonumber
E(|X-X'|) = 2 E(X)  + 2\int_{-\infty}^0[ 2F(x)  -F^2(x)]\, dx
-2\int_0^{+\infty}[ 1 -F(x)]^2\, dx,
\end{equation}
and, given two independent copies $Y$ and $Y'$ distributed according to $G$ 
\begin{equation}\nonumber
E(|Y-Y'|) = 2 E(Y)  + 2\int_{-\infty}^0[ 2G(x)  -G^2(x)]\, dx
-2\int_0^{+\infty}[ 1 -G(x)]^2\, dx.
\end{equation}
These identities clearly imply \fer{Cramer3}.

In one dimension, the Energy distance introduced by Szel\'ely \cite{Sze89, Sze03, SR} is nothing but twice the Cram\'er distance, as defined by \fer{Cramer3}. The interest in Energy distance is related to the fact that it can be expressed in any dimension, and that it allows for interesting generalizations. 

In what follows we refer to the Energy distance between two multidimensional random variables $\bX$ and $\bY$ with densities $f(\bx)$, and, respectively, $g(\bx)$, with $\bx \in \R^n$,  $n \ge 1$. The Energy distance of order $\alpha$, with $0<\alpha<2$ is given by
\begin{equations}\label{energy-a}
E_\alpha (f,g) &= 2\int_{\R^n\times\R^n}|\bx-\by|^\alpha f(\bx)g(\by) \, d\bx d\by\\
& -\int_{\R^n\times\R^n}|\bx-\by|^\alpha f(\bx)f(\by) \, d\bx d\by
 -\int_{\R^n\times\R^n}|\bx-\by|^\alpha g(\bx)g(\by) \, d\bx d\by\\
&= - \int_{\R^n\times\R^n}|\bx-\by|^\alpha [f(\bx)-g(\bx)][f(\by) -g(\by)] \, d\bx d\by.
\end{equations} 
As for the Cram\'er distance, the Energy distance of order $\alpha$ can be equivalently expressed in Fourier transform. In this case it holds \cite{ABGT2,SR}
\be\label{energy2}
E_\alpha (f,g) = c_{n,\alpha} \int_{\R^n} \frac{ |\ff(\bxi) -\gg(\bxi)|^2}{|\bxi|^{n+\alpha}}\, d\bxi,
\ee
where
\be\label{cn}
c_{n,\alpha} = \frac{\alpha 2^\alpha \Gamma\left( \frac{n+\alpha}2\right)}{2\pi^{n/2} \Gamma\left( \frac{2-\alpha}2\right)}.
\ee
The passage to Fourier transform is the easiest way to verify that right-hand side in \fer{energy-a} is nonnegative, so that the mixed integral dominates the sum of the last two. Further, expression \fer{energy2} allows one to relate in a precise way the energy distance of order $\alpha$ to homogeneous Sobolev spaces of fractional order.

\section{Convergence to equilibrium for one-dimensional Fokker--Planck equations}\label{sec:FP}

To better appreciate how  Cram\'er and Energy distances work in various problems related to kinetic equations, we will apply Cram\'er's distance to investigate the large time behavior of the solution $f(x,t)$ to the linear one-dimensional Fokker--Planck equation with constant diffusion $\sigma>0$, i.e.
\be\label{FoP}
\frac{\partial}{\partial t} f(x,t) = \frac{\partial}{\partial x} \left[ \sigma \frac{\partial}{\partial x}f(x,t) +  xf(x,t)\right],
\ee
when the initial value $f_0(x)$  is a probability density. Since equation \fer{FoP} is mass preserving, the solution remains a probability density for all subsequent times. Consequently, equation \fer{FoP} can be equivalently written in terms of the cumulative function $F(x,t)$, as defined in \fer{cumu}. Taking the integrals on both sides with respect to the spatial variable from $-\infty$ to $x$ we obtain that $F(x,t)$ satisfies 
\be\label{FoP-cum}
\frac{\partial}{\partial t} F(x,t) =  \sigma \frac{\partial f(x,t) }{\partial x}+  x\frac{\partial F(x,t)}{\partial x},
\ee
with initial value $F_0(x)$, the cumulative function of the initial   probability density $f_0(x)$.

If now $g(x,t)$ is another solution to \fer{FoP} starting from an initial probability density $g_0(x)$, in reason of the linearity, the difference between the two cumulative functions $F(x,t)$ and $G(x,t)$  solves the equation
\be\label{FoP-dif}
\frac{\partial}{\partial t}[F(x,t)-G(x,t)] =  \sigma \frac{\partial [f(x,t)-g(x,t)] }{\partial x}+  x\frac{\partial [F(x,t)-G(x,t)]}{\partial x}.
\ee
Therefore 
\begin{equations}\label{decay}
&\frac d{dt}\int_\R[F(x,t) -G(x,t)]^2\, dx = 2 \int_\R [F(x,t) -G(x,t)]\frac{\partial}{\partial t}[F(x,t)-G(x,t)]\, dx \\
&= 2\sigma \int_\R [F(x,t) -G(x,t)]\frac{\partial}{\partial x}[f(x,t)-g(x,t)]\, dx + \int_\R x\frac{\partial}{\partial x}[(F(x,t)-G(x,t))^2]\, dx \\
&= -2\sigma \int_\R\left[f(x,t)-g(x,t)\right]^2\, dx - \int_\R [F(x,t)-G(x,t)]^2\, dx.
\end{equations}
Note that the last line in \fer{decay} simply follows from integration by parts. 

These simple computations show that the Cram\'er distance between two solutions of the linear one-dimensional Fokker--Planck equation decays to zero at exponential rate, and that each term in it contributes to the decay. In particular, the linear drift is responsible of an exponential decay at rate $\exp\{-t\}$, while the diffusion term induces a decay that can be easily evaluated by resorting to some well-known inequality, in this case Nash inequality. We will not dwell further on these mathematical tools.

The important aspect to be outlined is that this result of decay continues to hold in the case in which $G(x,t)$ coincides with the equilibrium solution $F_\infty(x)$ of the Fokker--Planck equation, so that the Cram\'er distance provides (in few lines) a proof of (better than) exponential convergence towards equilibrium of the solution to equation \fer{FoP}. 

An interesting aspect of the previous computations is related to the fact that in the one-dimensional case it can be applied to a variety of Fokker--Planck type equations, characterized by a linear drift but different diffusion terms.  Consider for example the case of the Fokker--Planck equation
\be\label{poro}
\frac{\partial}{\partial t} f(x,t) = \frac{\partial}{\partial x} \left[  \frac{\partial}{\partial x}f^p(x,t) +  xf(x,t) \right],
\ee 
that, for $p>1$, is characterized by the Barenblatt equilibrium solution 
\be\label{Bar}
f_\infty(x) = \left(C^2- x^2\right)^{1/(p-1)},
\ee
and the constant $C$ is chosen so that the mass of $f_\infty$ is equal to $1$. Then the previous analysis applies with the main difference related to the fact that now the decay of Cram\'er's distance between two solutions is given by 
\begin{equations}\label{dec-p}
&\frac d{dt}\int_\R[F(x,t) -G(x,t)]^2\, dx \\ 
&= -2 \int_\R\left[f(x,t)-g(x,t)\right]\left[f^p(x,t)-g^p(x,t)\right]\, dx 
- \int_\R [F(x,t)-G(x,t)]^2\, dx.
\end{equations}
Since in \fer{dec-p} the contribution of the diffusion term is negative,  this  results in a (better than) exponential decay to zero of Cram\'er distance at  rate $\exp\{-t\}$. It is remarkable that a related improvement of the exponential decay for equation \fer{poro} has been observed before in \cite{CaTo1} by resorting to the decay of the relative Renyi entropies.

Maybe the most interesting results one can obtain through Cram\'er distance are related to one-dimensional Fokker--Planck equations of interest in socio-economic problems \cite{FPTT,NPT,PT13}. As far as the economic situation in a western society is concerned, a very popular mathematical model suitable to describe the evolution of the distribution of wealth $f(x,t)$, $x \in\R_+$, is  the  Fokker--Planck equation \cite{BM,CoPaTo05} 
 \be\label{FP2c}
 \frac{\partial f(x,t)}{\partial t} =  \frac \sigma{2}\frac{\partial^2 }
{\partial x^2}\left( x^2 f(x,t)\right) + \lambda \frac{\partial }{\partial x}\left(
(x-1) f(x,t)\right).
\ee
 The key feature of equation \fer{FP2c} is  that, for a given initial probability density $f_0(x)$ of mean value equal to $1$, and in presence of suitable boundary conditions at the point $x=0$, the solution preserves both mass and momentum and approaches in time a unique stationary solution of unitary mass \cite{TT1}.
This stationary state is given by the inverse Gamma distribution
 \be\label{equi2}
f_\infty(w) =\frac{(\mu-1)^\mu}{\Gamma(\mu)}\frac{\exp\left(-\frac{\mu-1}{w}\right)}{w^{1+\mu}},
 \ee
  where the positive constant $\mu >1$ is given by
\[
\mu = 1 + 2 \frac{\lambda}{\sigma}.
\]
As predicted by the Italian economist Vilfredo Pareto \cite{Par}, the equilibrium solution \fer{equi2} exhibits a power-law tail for large values of the wealth variable.

The rate of convergence to equilibrium of the solution to equation \fer{FP2c} was studied in \cite{ToTo} by expressing the linear kinetic equation in the Fourier space, that made possible to resort to a metric equivalent to the one-dimensional Energy distance of order $\alpha$, by suitably selecting $0<\alpha <2$.
As a matter of fact, if we apply to equation \fer{FP2c} the same procedure we applied to the classical Fokker--Planck equation \fer{FoP} we end up with the equation
\begin{equations}\label{dec-FP2c}
&\frac d{dt}\int_{\R_+}[F(x,t) -G(x,t)]^2\, dx \\ 
&= -\sigma \int_{\R_+}x^2 \left[f(x,t)-g(x,t)\right]^2\, dx 
-\lambda \int_{\R_+} [F(x,t)-G(x,t)]^2\, dx,
\end{equations}
which guarantees exponential convergence to equilibrium in Cram\'er's distance at rate $\exp\{-\lambda t\}$.

The last example we consider concerns with a mathematical model describing opinion formation \cite{Tos06}. There, the opinion variable is assumed to take values in  the bounded interval $ \CI=(-1, 1)$, the values $\pm 1$ denoting the extremal opinions. 
Among the various models introduced in \cite{Tos06}, one Fokker--Planck type equation has to be distinguished in view of its equilibrium configurations, which are represented by Beta-type probability densities supported in the interval $(-1, 1)$. This Fokker--Planck equation for the evolution of opinion density $f(x,t)$, with $|x| < 1$,  is given by 
 \be\label{op-FP}
 \frac{\partial f(x,t)}{\partial t} = \frac \lambda 2\frac{\partial^2 }{\partial x^2}\left((1-x^2)
 f(x,t)\right) + \frac{\partial }{\partial x}\left((x -m)f(x,t)\right).
 \ee
 In \fer{op-FP}, $\lambda/2 >0$ is the diffusion coefficient,  while $m$, with  $-1 <m<1$ is the mean opinion. 
Suitable boundary conditions at the boundary points $y = \pm 1$ then guarantee conservation of mass and momentum of the solution \cite{FPTT}.  

If the initial condition is a probability density, the steady state equals a probability density of Beta type, given by
 \be\label{beta}
f_{m,\lambda}(x)= C_{m,\lambda} (1-x)^{-1 + \frac{1-m}\lambda} (1+x)^{-1 + \frac{1+m}\lambda}.
 \ee
In \fer{beta} the constant $C_{m,\lambda}$ is such that the mass of $f_{m,\lambda}$
is equal to one. Since $-1 <m<1$, $f_{m,\lambda}$ is integrable on $\CI$.
Note that $f_{m,\lambda}$  is continuous on $\CI$, and as soon as $\lambda > 1+|m|$ tends to infinity as $y \to \pm 1$. 
Convergence to equilibrium in relative entropy for the solution to equation \fer{op-FP} has been studied in \cite{FPTT2}, where it was shown that entropy methods can  produce exponential convergence in $L_1$ towards equilibrium with an explicit rate, for some range of the parameters $\lambda$ and $m$ only. Owing to the weaker Cram\'er distance, we can easily fill up the gap, by showing that exponential convergence to equilibrium holds for all admissible values of the parameters. 
Indeed, the  procedure we applied to the classical Fokker--Planck equation \fer{FoP}  gives now
\begin{equations}\label{dec-op}
&\frac d{dt}\int_{-1}^1 [F(x,t) -G(x,t)]^2\, dx \\ 
&= -\frac\lambda{2} \int_{-1}^1(1-x^2) \left[f(x,t)-g(x,t)\right]^2\, dx 
- \int_{-1}^1 [F(x,t)-G(x,t)]^2\, dx,
\end{equations}
thus showing that, at difference with the results in \cite{FPTT2},  the (better than) exponential convergence to equilibrium in Cram\'er distance holds independently of the values of the parameters. 

Clearly, in view of identity \fer{Cramer3}, the same results can be obtained by resorting to the Energy distance of order $\alpha =1$. It is important to remark that the choice of an Energy distance of order $\alpha \not=1$ leads to more intricate computations. An example of this can be found in \cite{TT1}, where the results of exponential convergence in a certain range of the parameter $\alpha$ have been obtained by resorting to the Fourier expression of the distance. We will be back on this aspect in the next Section.

\section{The linear Fokker--Planck equation in $\R^n$}\label{sec:RN}

In this Section we investigate the rate of decay of the Energy distance of order $\alpha$ along the solution to the linear $n$-dimensional Fokker--Planck equation, $n>1$
\be\label{FP}
\frac{\partial}{\partial t} f(\bx,t) = \nabla\cdot \left[ \nabla f(\bx,t) + \bx f(\bx,t) \right].
\ee
For the sake of simplicity, we will study separately the decay of the Energy distance on the linear drift equation
\be\label{drift}
\frac{\partial}{\partial t} f(\bx,t) = \nabla\cdot \left[ \bx f(\bx,t)) \right]
\ee
and the linear diffusion equation
\be\label{heat}
\frac{\partial}{\partial t} f(\bx,t) = \Delta f(\bx,t).
\ee
Let $0<\alpha<2$. Thanks to definition \fer{energy-a} we have 
\begin{equations}\label{prop1}
\frac d{dt} E_\alpha (f(t),g(t)) 
&= - \frac d{dt} \int_{\R^n\times\R^n}|\bx-\by|^\alpha [f(\bx,t)-g(\bx,t)][f(\by,t) -g(\by,t)] \, d\bx d\by \\
&=- \int_{\R^n\times\R^n}|\bx-\by|^\alpha [f(\by,t)-g(\by,t)]\frac\partial{\partial t}[f(\bx,t) -g(\bx,t)] \, d\bx d\by\\
&- \int_{\R^n\times\R^n}|\bx-\by|^\alpha [f(\bx,t)-g(\bx,t)]\frac\partial{\partial t}[f(\by,t) -g(\by,t)] \, d\bx d\by. 
\end{equations} 
If we apply identity \fer{prop1} to the drift equation \fer{drift} and integrate by parts,  we obtain 
\begin{equations}\label{drift1}
& \int_{\R^n\times\R^n}|\bx-\by|^\alpha [f(\by,t)-g(\by,t)]\frac\partial{\partial t}[f(\bx,t) -g(\bx,t)] \, d\bx d\by\\
&=\int_{\R^n\times\R^n}|\bx-\by|^\alpha [f(\by,t)-g(\by,t)]\nabla_\bx\cdot[\bx( f(\bx,t) -g(\bx,t))] \, d\bx d\by \\
&= -\int_{\R^n\times\R^n}\bx\cdot \nabla_\bx |\bx-\by|^\alpha [f(\by,t)-g(\by,t)][ f(\bx,t) -g(\bx,t)] \, d\bx d\by \\
&= -\alpha \int_{\R^n\times\R^n}|\bx-\by|^{\alpha-2}\sum_{k=1}^n x_k(x_k-y_k) [f(\by,t)-g(\by,t)][ f(\bx,t) -g(\bx,t)] \, d\bx d\by.
\end{equations} 
Likewise
\begin{equations}\label{drift2}
&\int_{\R^n\times\R^n}|\bx-\by|^\alpha [f(\bx,t)-g(\bx,t)]\frac\partial{\partial t}[f(\by,t) -g(\by,t)] \, d\bx d\by\\
&= -\alpha \int_{\R^n\times\R^n}|\bx-\by|^{\alpha-2}\sum_{k=1}^n y_k(y_k-x_k) [f(\by,t)-g(\by,t)][ f(\bx,t) -g(\bx,t)] \, d\bx d\by.
\end{equations}
Therefore, taking the sum of \fer{drift1} and\fer{drift2} we obtain that the drift equation \fer{drift} satisfies
\be\label{E-dec}
\frac d{dt} E_\alpha (f(t),g(t)) = -\alpha E_\alpha (f(t),g(t)).
\ee
Hence, along the drift equation, the Energy distance of order $\alpha$ between two solutions decays exponentially at a rate proportional to its order.

Let us now apply identity \fer{prop1} to the linear diffusion equation \fer{heat} and integrate twice by parts. We obtain
\begin{equations}\label{heat1}
&\int_{\R^n\times\R^n}|\bx-\by|^\alpha [f(\by,t)-g(\by,t)]\Delta[ f(\bx,t) -g(\bx,t))] \, d\bx d\by \\
&= \int_{\R^n\times\R^n}\left(\Delta_\bx |\bx-\by|^\alpha\right) [f(\by,t)-g(\by,t)][ f(\bx,t) -g(\bx,t)] \, d\bx d\by \\
&= \alpha(n-2+\alpha) \int_{\R^n\times\R^n}\frac 1{|\bx-\by|^{2-\alpha}} [f(\by,t)-g(\by,t)][ f(\bx,t) -g(\bx,t)] \, d\bx d\by.
\end{equations} 
Likewise
\begin{equations}\label{heat2}
&\int_{\R^n\times\R^n}|\bx-\by|^\alpha [f(\bx,t)-g(\bx,t)]\Delta[ f(\by,t) -g(\by,t))] \, d\bx d\by \\
&= \alpha(n-2+\alpha) \int_{\R^n\times\R^n}\frac 1{|\bx-\by|^{2-\alpha}} [f(\by,t)-g(\by,t)][ f(\bx,t) -g(\bx,t)] \, d\bx d\by.
\end{equations}

For $0< \alpha <2$, let us denote
    \begin{equation}
        \label{wlambda}
        W_\alpha(\bx) = |\bx|^{-(2-\alpha)}.
    \end{equation}
    As shown in Theorem 4.1 from \cite{Stein}, since for $n \ge 2$ one has $2-\alpha <n$,  the Fourier transform of $W_\alpha$ equals
    \begin{equation}
        \label{fW}
        \widehat W_\alpha(\bxi) = \int_{\R^n} |\bx|^{-(2-\alpha)} e^{-i\bxi\cdot\bx} \, d\bx = \pi^{n/2}2^{n-2 + \alpha} \frac{\Gamma\left(\frac{n-2 + \alpha}2\right)}{\Gamma\left(\frac{2-\alpha}2\right)}\, \frac 1{|\bxi|^{n-2 +\alpha}} 
    \end{equation}
    where $\Gamma(\cdot)$ is the gamma function. This implies \cite{ABGT2}
    \begin{equations}\label{energy-}
&0 \le E_{-(2-\alpha)}(f(t),g(t)) = \int_{\R^n\times\R^n}\frac 1{|\bx-\by|^{2-\alpha}} [f(\by,t)-g(\by,t)][ f(\bx,t) -g(\bx,t)] \, d\bx d\by  \\
& = d_{n,\alpha}\int_{\R^n}\frac{ |\ff(\bxi,t)-\gg(\bxi,t)|^2} {|\bxi|^{n-2+ \alpha}}\, d\bxi,
\end{equations} 
where
\be\label{dn}
d_{n,\alpha} = \frac{2^\alpha}{2^2\pi^{n/2}} \frac{\Gamma\left(\frac{n-2 + \alpha}2\right)}{\Gamma\left(\frac{2-\alpha}2\right)}
\ee
The Energy distance of negative order $2-\alpha$, with $0<\alpha<2$, defined in \fer{energy-} has been introduced and studied in \cite{ABGT2}. It is interesting to remark that, at difference with the Energy distance of positive order, the mixed integral is now dominated by the sum of the other two. Indeed, expression \fer{energy-} can be equivalently expressed as 
\begin{equations}\label{energy-2-a}
E_{-(2-\alpha)} (f,g) &= 
\int_{\R^n\times\R^n}\frac 1{|\bx-\by|^{2-\alpha}} f(\bx)f(\by) \, d\bx d\by \\
 &+\int_{\R^n\times\R^n}\frac 1{|\bx-\by|^{2-\alpha}} g(\bx)g(\by) \, d\bx d\by - 2\int_{\R^n\times\R^n}\frac 1{|\bx-\by|^{2-\alpha}} f(\bx)g(\by) \, d\bx d\by.
\end{equations} 

Joining the two results, if $f(\bx,t)$ solves the Fokker--Planck equation \fer{FP} and $f_\infty(\bx)$  denotes the steady solution Gaussian density we proved that
\begin{equations}\label{dec-n}
&\frac d{dt} E_\alpha (f(t),f_\infty) 
&=-2 \,\alpha(n-2+\alpha)E_{-(2-\alpha)}(f(t),f_\infty) -\alpha  E_\alpha (f(t),f_\infty).
\end{equations}

\section{The Energy distance and the linear diffusion equation}\label{sec:heat}

It is interesting to remark that, for any given $0<\alpha<2$, and $n\ge 2$, given two probability densities $f(\bx)$ and $g(\bx)$, $\bx \in \R^n$, the Energy distance $E_{-(2-\alpha)}(f,g)$, defined as in \fer{energy-},  controls the Energy distance $E_\alpha(f,g)$ defined in \fer{energy2}. An explicit bound can be found as follows.  Let $R>0$. Then
\be\label{e1}
\int_{|\bxi| \le R}\frac{ |\ff(\bxi)-\gg(\bxi)|^2} {|\bxi|^{n+ \alpha}}\, d\bxi \le \sup_{\bxi}\frac{ |\ff(\bxi)-\gg(\bxi)|^2} {|\bxi|^{2}}\int_{|\bxi| \le R}\frac{1} {|\bxi|^{n+ \alpha-2}}\, d\bxi .
\ee
Thanks to Proposition $2.6$ from \cite{CaTo}, we know that, for any pair of probability densities $f(\bx)$ and $g(\bx)$, the Fourier--based distance
\[
d_1(f,g) = \sup_{\bxi \in \R^n}\frac{ |\ff(\bxi)-\gg(\bxi)|} {|\bxi|} 
\]
is bounded. Since $\alpha <2$, by solving the integral in \fer{e1} we obtain
\[
\int_{|\bxi| \le R}\frac{1} {|\bxi|^{n+ \alpha-2}}\, d\bxi = \frac{2\pi^{n/2}}{\Gamma\left(\frac n2 \right)}\int_0^R \rho^{1-\alpha}\,d\rho = \frac{2\pi^{n/2}}{\Gamma\left(\frac n2 \right)}\frac{R^{2-\alpha}}{2-\alpha}.
\]
Therefore, if $c_{n,\alpha}$ is the constant defined in \fer{cn}, 
\be\label{e2}
c_{n,\alpha} \int_{|\bxi| \le R}\frac{ |\ff(\bxi)-\gg(\bxi)|^2} {|\bxi|^{n+ \alpha}}\, d\bxi \le c_{n,\alpha}\frac{2\pi^{n/2}}{\Gamma\left(\frac n2 \right)} d_1(f,g)^2 \frac{R^{2-\alpha}}{2-\alpha}.
\ee
Further we have
\[
\int_{|\bxi| > R}\frac{ |\ff(\bxi)-\gg(\bxi)|^2} {|\bxi|^{n+ \alpha}}\, d\bxi =\int_{|\bxi| > R}\frac 1{|\bxi|^2} \frac{ |\ff(\bxi)-\gg(\bxi)|^2} {|\bxi|^{n+ \alpha-2 }}\, d\bxi \le \frac 1{R^2} \int_{|\bxi| > R}\frac{ |\ff(\bxi)-\gg(\bxi)|^2} {|\bxi|^{n+ \alpha-2}}\, d\bxi.
\]
Consequently
\be\label{e3}
c_{n,\alpha} \int_{|\bxi| > R}\frac{ |\ff(\bxi)-\gg(\bxi)|^2} {|\bxi|^{n+ \alpha}}\, d\bxi \le \frac{c_{n,\alpha} }{d_{n,\alpha} }E_{-(2-\alpha)}(f,g)\frac 1{R^2}.
\ee
Joining \fer{e2} and \fer{e3} we obtain
\be\label{e4}
E_\alpha(f,g) \le A_{n,\alpha} d_1(f,g)^2 \frac{R^{2-\alpha}}{2-\alpha} + B_{n,\alpha} E_{-(2-\alpha)}(f,g) \frac 1{R^2}, 
\ee
with obvious meaning of the constants $A_{n,\alpha}$ and $B_{n.\alpha}$. Optimizing over $R$ in \fer{e4} we finally obtain
\be\label{e-fin}
E_\alpha(f,g) \le D_{n,\alpha} d_1(f,g)^{\frac 4{4-\alpha}} E_{-(2-\alpha)}(f,g)^{\frac{2-\alpha}{4-\alpha}},
\ee
where
\be\label{Dn}
 D_{n,\alpha} = A_{n,\alpha}^{\frac 2{4-\alpha}} B_{n,\alpha}^{\frac{2-\alpha}{4-\alpha}}\left[ \frac{2^{\frac{2-\alpha}{4-\alpha}}}{2-\alpha}+ \left(\frac 12\right)^{\frac 2{4-\alpha}} \right].
\ee
%
Let $f(t)$ and $g(t)$ be two solutions to equation \fer{heat} with initial data $f_0(\bx)$ and $g_0(\bx)$, both being probability densities on $\R^n$.
By examining the Fourier transform formulation of equation \fer{heat}, it is straightforward to verify that the distance $d_1(f(t),g(t))$ is monotone non-increasing in time.
Applying inequality \fer{e-fin} to the linear diffusion equation \fer{heat}, we then obtain from \fer{dec-n} that
\begin{equation}\label{decay-n}
\frac{d}{dt} E_\alpha(f(t),g(t))
= -2\,\alpha(n-2+\alpha)\,E_{-(2-\alpha)}(f(t),g(t))
\le -C\,E_\alpha(f(t),g(t))^{\frac{4-\alpha}{2-\alpha}}.
\end{equation}
where
\be\label{cc} 
C=  C(\alpha,n, d_1(f_0,g_0))= 2\,\alpha(n-2+\alpha) \left[ D_{n,\alpha} d_1(f_0,g_0)^{\frac 4{4-\alpha}} \right]^{- \frac{4-\alpha}{2-\alpha}}.
\ee
Integrating \fer{decay-n} we obtain
\be\label{new}
E_\alpha (f(t),g(t))  \le \left( \frac 1{E_\alpha (f_0,g_0)^{(2-\alpha)/2}} + C\frac 2{2-\alpha} t \right)^{- {\frac{2-\alpha}{2}}}.
\ee
Hence the Energy distance of order $\alpha$  between two different solutions decays to zero at a rate $t^{- {\frac{2-\alpha}{2}}}$, a polynomial rate decreasing with respect to $\alpha$. Clearly, the same result remains valid if we choose as $g(\bx,t)$ the fundamental solution of the heat equation, namely the Gaussian density of zero mean and variance $2t$.

\section{Conclusions}\label{sec:conclusions}

In this paper,  we studied the rate of convergence to equilibrium of the solutions to  Fokker--Planck type equations  with linear drift by means of Cram\'er and Energy distances.   In all cases in which the Fokker--Planck type equations can be treated through these distances, it is shown that the rate of decay is improved with respect to known results which are based on the decay of relative entropy. It is interesting to remark that the study of the decay to equilibrium by means of these distances is not restricted to Fokker--Planck type equations, but it can be usefully applied to systems of Fokker--Planck type equations, as recently done in \cite{MTZ}, or to the linear diffusion equation, as briefly shown in Section \ref{sec:heat}. The present results suggest that both Cram\'er and Energy distances could be of high interest in the study of both theoretical and numerical problems related to Fokker--Planck type equations.


\section{acknowledgement}
This work has been written within the activities of GNFM (Gruppo Nazionale per la Fisica Matematica) of INdAM (Istituto Nazionale di Alta Matematica), Italy.
 The results contained in the present paper have been partially presented in WASCOM 2025. The authors state that there is no conflict of interest.

\end{document}